\documentclass[english,russian]{amsart}
\usepackage{amssymb,amsmath}

\newtheorem{thm}{Theorem}%[section]
\newtheorem{crl}[thm]{Corollary}

\newtheorem{lm}[thm]{Lemma}

\newcommand{\der}{\partial}

\newcommand{\sudda}[1]{}

\newcommand{\limsim}{\limits_{\sim\!\sim\!\sim\!\sim\!\sim\!\sim\sim\!\sim\!\sim\!\sim\!\sim\!\sim}}

\newcommand{\limsimeq}{\limits_{\simeq\!\simeq\!\simeq\!\simeq\!\simeq\!\simeq
\simeq\!\simeq\!\simeq\!\simeq\!\simeq\!\simeq}}

\begin{document}

\title{Euclid and Frobenius structures on Weyl algebra}

\author{A.S. Dzhumadil'daev}

\address
{Kazakh-British University, Tole bi, 59, Almaty,  050000,
Kazakhstan} \email{dzhuma@hotmail.com} 

\maketitle

\begin{abstract}  
A non-degenerate associative bilinear form and a positive definite symmetric bilinear form on Weyl algebra are constructed. 

 \end{abstract}
 
 Weyl algebra  $A_n$ is defined as an associative algebra over a field $K$ of characteristic $0$ 
 generated by $2n$ generators $x_1,\ldots,x_n,y_1,\ldots,y_n$ with defining relations $[x_i,y_j]=\delta_{i,j},$ where $\delta_{i,j}$ is Kronecker symbol, $\delta_{i,i}=1$ and $\delta_{i,j}=0,$ if $i\ne j.$ In our paper we wiil use interpretation of Weyl algebra in terms of differential operators. Details see \cite{Bjork}. 
  
 Let ${\bf Z}$ be integers group, ${\bf Z}_0$ its subset of non-negative integers and ${\bf Z}_0^n=
 \{\alpha=(\alpha_1,\ldots,\alpha_n) | \alpha_i\in {\bf Z}_0\}.$   Let $U_n=K[x_1,\ldots,x_n]$ be an algebra of polynomials with $n$ variables. It has base collected by $x$-monoms  $x^\alpha=x_1^{\alpha_1}\cdots x_n^{\alpha_n},$ where $\alpha\in {\bf Z}_0^n.$ Let 
$\der_i:U_n\rightarrow U_n, 1\le i\le n,$ be a derivation of $U_n$ defined by  
 $$\der_i x^{\alpha}=\alpha_i x^{\alpha-\epsilon_i},$$
 where $\epsilon_i=(0,\ldots,0,{1},0,\ldots,0)\in {\bf Z}_0^n$ (all components, except $i$-th are $0$).  Let  $\der^{\alpha}=\der_1^{\alpha_1}\cdots\der_n^{\alpha_n}$ be $\der$-monom.   Call a product of $x$-monom and $\der$-monom  $x^\alpha\der^\beta,$ where $\alpha,\beta\in {\bf Z}_0^n,$ as a {\it  Weyl monom} or simply monom. Let $A_n$ be a linear  span of Weyl monoms. Note that a linear combination of Weyl monoms can be interpretered as a differential operator on $U_n.$ Therefore, we can endow $A_n$ by a structure of associative algebra under composition operatotion defined of Weyl monoms by 
 $$x^{\alpha}\der^{\beta}\circ  x^{\alpha'}\der^{\beta'}=\sum_{\gamma\in {\bf Z}_0^n}\gamma! {\beta\choose \gamma}{\alpha'\choose \gamma}x^{\alpha+\alpha'-\gamma}\der^{\beta+\beta'-\gamma}.$$
 Here we set
 $$\alpha!=\alpha_1!\cdots \alpha_n!,\quad {\alpha\choose \beta}={\alpha_1\choose \beta_1}\cdots
 {\alpha_n\choose \beta_n}.$$

For a Weyl monom  $X=x^{\alpha}\der^{\beta}\in A_n$ say that $X$ has {\it  multi-weight} $w(X)=\alpha-\beta\in {\bf Z}^n$ and {\it weight} $l(X)=\sum_{i=1}^n \alpha_i-\beta_i.$ Let $A_{n,\omega}$ be a subspace of $A_n$ generated by monoms of weight $\omega.$ 
Say that $X\in A_n$ is  {\it a differential operator of multi-weight} $\omega,$ and write $w(X)=\omega,$  if $X\in A_{n,\omega}.$ 
 Similarly, let $A_{n,l}$ be subspace of $A_n$ generated by  monoms of weight $l.$ Say that $X\in A_n$ is  {\it a differential operator of weight} $l,$ and write $l(X)=l,$  if $X\in A_{n,l}.$ Note that for $n=1$ notions multi-weight and weight are coincide. 
 
{ \bf Proposition.} {\it  If $X$ and $X'$ are differential operators of multi-weight $\omega$ and $\omega',$ then $X\circ X'$ is a differential operator of multi-weight  $\omega+\omega'.$ If $X$ and $X'$ are differential operators of weight $l$ and $l',$ then $X\circ X'$ is a differential operator of weight $l+l'.$ }

Therefore we can endow $A_n$ by multi-weight ${\bf Z}^n$-grading  and by weight ${\bf Z}$-grading,
$$A_n=\oplus_{\omega\in {\bf Z}^n} A_{n,\omega},\qquad A_{n,\omega}\circ A_{n,\omega'}\subseteq A_{n,\omega+\omega'},$$
 $$A_n=\oplus_{l\in {\bf Z}} A_{n,l},\qquad A_{n,l}\circ A_{n,l'}\subseteq A_{n,l+l'}.$$
 Note that homogeneous components $A_{n,\omega}$ and $A_{n,l}$ are infinite-dimensional. Note also that $A_{n,0}$ is a subalgebra of $A_n.$ For $X\in A_n$ denote by $X_\omega$ its projection on $A_{n,\omega}.$  
% Let $X\mapsto X_{\omega}$ be a projection map  $A_n\rightarrow A_{n,\omega}.$ 
 
 For a monom $X=x^{\alpha}\der^{\beta}$ denote by $X^{(i)}$ its part on 
 $$A_1^{(i)}=\{x_i^{\alpha_i}\der_i^{\beta_i} |\alpha_i,\beta_i\in {\bf Z}_0\}.$$
So, 
$$X=x^{\alpha}\der^{\beta}\Rightarrow X^{(i)}=x_i^{\alpha_i}\der_i^{\beta_i}.$$
Say that  $X$ is {\it decomposable}, if it can be presented as a product 
$X=\prod_{i=1}^n X^{(i)},$ where $X^{(i)}\in A_1^{(i)}.$ Since any $X\in A_n$ can be presented as a linear combination of decomposable elements, a linear map 
$$A_1^{(1)}\otimes \cdots \otimes A_1^{(n)}\rightarrow A_n$$
defined on monoms by 
$$ X^{(1)}\otimes \cdots \otimes X^{(n)}\mapsto X^{(1)}\cdots X^{(n)}$$
induces an isomorphism of algebras.

A linear map $X\rightarrow \overline X$ given by 
$$\overline{x^{\alpha}\der^\beta} =x^\beta \der^\alpha,$$
 generates anti-authomorphsim of Weyl algebra,
 $$\overline{X\circ Y}=\overline Y\circ \overline X.$$
 
 An operator $X\in A_n$ is called {\it self-adjoint} if  $\overline X=X$  and {\it skew-adjoint,}  if 
 $\overline X=-X.$ Let $A_n^{(+)}$ and $A_n^{(-)}$ are subspaces of $A_n$ generated by self-adjoint and skew-adjoint operators.

 Let us define a linear map $T: A_n\rightarrow K$ by 
 $$T(x^{\alpha}\der^{\beta})=\left\{ \begin{array}{cc}
 \alpha!& \mbox{ if $\alpha=\beta$}\\
 0&\mbox{ otherwise}\end{array}\right.
 $$
 Let $(\;\, ,\; )$ be a bilinear form on $A_n$ defined by 
 $$(X,Y)=T(X\circ Y).$$

 Recall that an associative algebra is called {\it Frobenius} if it has a symmetric non-degenerate  associative bilinear form and  {\it Euclid} if it has a symmetric positive definite bilinear form. 
 
 Our aim is to establish the following results 
 
 \begin{thm} \label{Frobenius}
 The bilinear form $(\;\, , \;): A_n\times A_n \rightarrow K$ is 
 \begin{itemize}
 \item symmetric in the following sences: 
 \begin{equation}\label{51}
 (X,Y)=2^{l(Y)} (Y,X),
 \end{equation}
\begin{equation}\label{52}
\sqrt{2}^{l(X)}(X,Y)=\sqrt{2}^{l(Y)}(Y,X),
\end{equation}
 for any homogeneous $X, Y\in A_n,$ and 
 \begin{equation}\label{53}
 (\overline X, \overline Y)=(Y,X),
 \end{equation}
 for any $X, Y\in A_n.$
 \item assosiative
 $$(X\circ Y,Z)=(X,Y\circ Z),$$
 for any $X,Y,Z\in A_n.$
 \item non-degenerate,
 $$(X,Y)=0, \quad  \forall Y\in A_n\Rightarrow X=0.$$
 \end{itemize}
 In particular,   $A_{n,0}$, the  algebra of differential operators of weight $0,$ is Frobenius. 
 \end{thm}

 \begin{thm} \label{Euclid}
 The bilinear form $\langle \;\, , \;\rangle: A_n\times A_n \rightarrow K$ 
 given by 
 $$\langle X,Y\rangle=\sqrt{2}^{\,l(X) }(X,\overline{Y})$$
 has the following properties:
 \begin{itemize}
 \item symmetric 
 \begin{equation}\label{61}
 \langle X,Y\rangle =\langle Y, X\rangle, \qquad \forall X,Y\in A_n,
 \end{equation}
  \item one more symmetry  
  \begin{equation}\label{62}
 \langle X,Y\rangle =\langle \overline X,\overline Y\rangle,  \qquad \forall X,Y\in A_n,
 \end{equation}
 \item positive definite,
 $$\langle X,X\rangle \ge 0, \quad \forall X\in A_n,$$
 $$\langle X,X\rangle=0 \Rightarrow X=0.$$
 \item invariant in the following sence
 \begin{equation}\label{63}
 \langle a\circ X,Y\rangle =\sqrt{2}^{\,-l(a)}\langle X, \bar a\circ Y\rangle, 
 \end{equation}
 \begin{equation}\label{64}
 \langle X\circ a, Y\rangle=
 \sqrt{2}^{\, l(a)}\langle X,  Y\circ \bar a\rangle,
 \end{equation}
 for any $X,Y\in A_n$ and for any (weight)-homogenous $a\in A_n.$
 \item (multi)-weight subspaces are mutually orthogonal,
 $$X\in A_{n,\omega}, \; X'\in A_{n,\omega'}, \;  \omega\ne \omega'\Rightarrow \langle X,X'\rangle =0.$$
 In particular, 
 $$\langle X,Y\rangle=\sum_{\omega\in {\bf Z}^n} \langle X_{\omega},Y_{\omega}\rangle,$$
 and
 $$|X|^2=\sum_{\omega\in {\bf Z}^n} |X_{\omega}|^2,$$
 where $|X|=\sqrt{\langle X,X\rangle}$ is a length of $X.$ 
\item  for decomposable elements $X=X^{(1)}\cdots X^{(n)}$ and $Y=Y^{(1)}\cdots Y^{(n)},$ 
 $$\langle X,Y\rangle =\prod_{i=1}^n \langle X^{(i)},Y^{(i)}\rangle.$$
  \end{itemize}
  \end{thm}
  
  \begin{crl} The spaces of self-adjoint  and skew-adjoint differential  operators  $A_{n}^{(+)}$  and $A_n^{(-)}$ are  Euclid.  The subspaces $A_n^{(+)}$ and $A_n^{(-)}$ are mutually orthogonal,
 $$<X,Y>=0 ,\qquad \forall X\in A_n^{\pm}, \; \forall Y\in A_n^{\mp}.$$
  \end{crl}
 
 \begin{crl} For any $X,Y\in A_n$ and for homogeneous $a\in A_n,$ 
 $$\langle [a,X],Y\rangle=\langle X,\sqrt{2}^{\; -l(a)} \overline{a}\circ Y-\sqrt{2}^{\; l(a)}Y\circ \overline{a}\rangle.$$
 In particular, for any $a\in A_{n,0},$ $X,Y\in A_n,$ 
 $$\langle [a,X],Y\rangle=\langle X,[\, \overline{a},Y]\rangle.$$
 \end{crl}
    
 {\bf Remark.}  A norm induced by scalar product $|X|=\sqrt{\langle X,X\rangle}$ satisfies the triangle inequality
 $$|X+Y|\le |X|+|Y|,\qquad \forall X,Y\in A_n.$$
 By  normed algebra one can understand an algebra  $A=(A,\circ)$  with a norm $X\mapsto |X|$ that  satisfies the inequality  $$|X\circ Y|\le |X| |Y|,$$ for any $X,Y\in A.$ Weyl algebra $A_n$ under norm $|X|=\sqrt{\langle X,X\rangle}$ is not  normed in this sence. For example, $|x|=0, |\der|=0,$ but $|x\circ \der|=|x\der|=3\ge 0=|x||\der|.$  We think that 
 $$|X\circ Y|\ge |X| |Y|,$$ for any $X,Y\in A_n.$
 
 {\bf Example.} For any non-negative integer $k$ the scalar product 
 $$\langle(x\der)^i, (x\der)^{k-i}\rangle =\sum_{j=0}^k s_{k,j}j!$$
does not depend on $i.$ It is equal to the number of ordered partitions of a set with $k$ elements. 
 Here $s_{k,j}$ are Stirling numbers of second kind.  In particular,
 $$|(x\der)^k|=\sqrt{\sum_{j=0}^{2k} s_{2k,j}j!}.$$

 \section{Frobenius structure on Wel algebra}
 
 In this section we prove Theorem \ref{Frobenius}. To do that we need several lemmas. 
 
 \begin{lm}\label{1}
 $$\sum_{i,a,b} {n\choose b}{b\choose i}{n-i \choose a}(x+1)^a x^i y^az^b=(1+(1+x)(y+z+yz))^n.$$
 \end{lm}
 
 {\bf Proof.} Let us denote left hand of our relation by $G_n$ and right hand by $H_n.$ 
 We will prove by induction ion $n$  that $G_n=H_n.$
For $n=0$ is nothing to prove.

Suppose that $G_n=H_n$ for $n.$  Then
 $$H_{n+1}=$$ 
 
 $$H_n(1+(x+1)(y+z+yz))=G_n(1+(x+1)(y+z+yz))=$$
 
  $$\sum_{i,a,b} {n\choose b}{b\choose i}{n-i \choose a}(x+1)^a x^i y^az^b(1+(x+1)(y+z+yz))=$$

  \bigskip
    
  $$\sum_{i,a,b} {n\choose b}{b\choose i}{n-i \choose a}(x+1)^a x^i y^az^b+
  {n\choose b}{b\choose i}{n-i \choose a}(x+1)^{a+1} x^i y^{a+1} z^b+
$$
$$
 {n\choose b}{b\choose i}{n-i \choose a}(x+1)^{a+1} x^i y^a z^{b+1}+
 {n\choose b}{b\choose i}{n-i \choose a}(x+1)^{a+1} x^i y^{a+1}z^{b+1}  =$$
  
  \bigskip (in third summand take $(x+1)^{a+1}=(x+1)(x+1)^a$)

  $$\sum_{i,a,b} {n\choose b}{b\choose i}{n-i \choose a}(x+1)^a x^i y^az^b+
  {n\choose b}{b\choose i}{n-i \choose a}(x+1)^{a+1} x^i y^{a+1} z^b+
$$
$$
 {n\choose b}{b\choose i}{n-i \choose a}(x+1)^{a} x^i y^a z^{b+1}+
 {n\choose b}{b\choose i}{n-i \choose a}(x+1)^{a} x^{i+1} y^a z^{b+1}+$$
 $$
 {n\choose b}{b\choose i}{n-i \choose a}(x+1)^{a+1} x^i y^{a+1}z^{b+1}  =$$

  \bigskip
    
  $$\sum_{i,a,b} {n\choose b}{b\choose i}{n-i \choose a}(x+1)^a x^i y^az^b+
  {n\choose b}{b\choose i}{n-i \choose a-1}(x+1)^{a} x^i y^{a} z^b+
$$
$$
 {n\choose b-1}{b-1\choose i}{n-i \choose a}(x+1)^{a} x^i y^a z^{b}+
 {n\choose b-1}{b-1\choose i-1}{n-i+1 \choose a}(x+1)^{a} x^{i} y^a z^{b}+$$
 $$
 {n\choose b-1}{b-1\choose i}{n-i \choose a-1}(x+1)^{a} x^i y^{a}z^{b}  =$$

 \bigskip

  $$\sum_{i,a,b} \{\mathop{{n\choose b}{b\choose i}{n-i \choose a}}\limsim+
 \mathop{ {n\choose b}{b\choose i}{n-i \choose a-1}}\limsim+
$$
$$
\mathop{ {n\choose b-1}{b-1\choose i}{n-i \choose a}}\limsimeq+
 {n\choose b-1}{b-1\choose i-1}{n-i+1 \choose a}+$$
 $$\mathop{{n\choose b-1}{b-1\choose i}{n-i \choose a-1}}\limsimeq \} (x+1)^{a} x^i y^{a}z^{b}  =$$

 \bigskip

  $$\sum_{i,a,b} \{{n\choose b}{b\choose i}{n-i+1 \choose a}+
 {n\choose b-1}{b-1\choose i}{n-i+1 \choose a}+
 {n\choose b-1}{b-1\choose i-1}{n-i+1 \choose a} \} (x+1)^{a} x^i y^{a}z^{b}  =$$

 \bigskip

  $$\sum_{i,a,b} \{{n\choose b}{b\choose i}+
 {n\choose b-1}{b-1\choose i}+
 {n\choose b-1}{b-1\choose i-1}\}{n-i+1 \choose a} ( x+1)^{a} x^i y^{a}z^{b}  =$$

 \bigskip

  $$\sum_{i,a,b} \{{n\choose b}{b\choose i}+
 {n\choose b-1}{b\choose i}\}{n-i+1 \choose a} \}( x+1)^{a} x^i y^{a}z^{b}  =$$

 \bigskip

  $$\sum_{i,a,b} \{{n+1\choose b}{b\choose i}{n-i+1 \choose a} ( x+1)^{a} x^i y^{a}z^{b}  =$$

$$G_{n+1}.$$
So, the condition $G_n=H_n$ implies that $G_{n+1}=H_{n+1}.$

Lemma is proved. 

\begin{crl} \label{2} For any $a,b\in {\bf Z}_0,$
$$\sum_{i} 2^a{n\choose b}{b\choose i}{n-i\choose a}=
\sum_{i} 2^b{n\choose a}{a\choose i}{n-i\choose b}.$$
\end{crl}

{\bf Proof}. Note that the coefficient at $ y^a z^b$ of $G_n(x,y,z)$ is 
$\sum_{i} (x+1)^a{n\choose b}{b\choose i}{n-i\choose a}x^i.$ By Lemma \ref{1} 
it is equal to the coefficient at $y^b z^a$ of $G_n(x,y,z).$ So,
$$\sum_{i} (x+1)^a{n\choose b}{b\choose i}{n-i\choose a}x^i=
\sum_{i} (x+1)^b{n\choose a}{a\choose i}{n-i\choose b}x^i.$$
This relation  for $x=1$ gives us our statement.

\begin{lm} \label{3}For any $\alpha,\beta,\theta\in{\bf Z}_0^n,$
$$\sum_{\gamma\in {\bf Z}^n_0} 2^{||\alpha||}\frac{{\theta-\gamma\choose \alpha}{\beta\choose 
\gamma}}{{\theta\choose \alpha}}=
\sum_{\gamma\in {\bf Z}^n_0} 2^{||\beta||}\frac{{\theta-\gamma\choose \beta}{\alpha\choose 
\gamma}}{{\theta\choose \beta}},$$
where $||\alpha||=\sum_{i=1}^n \alpha_i$ and $||\beta||=\sum_{i=1}^n \beta_i.$
\end{lm}

{\bf Proof.} Let $\alpha=(\alpha_1,\ldots,\alpha_n), \beta=(\beta_1,
\ldots,\beta_n), \gamma=(\gamma_1,\ldots,\gamma_n),\theta=(\theta_1,\ldots,\theta_n).$ Then 

$$\sum_{\gamma\in {\bf Z}^n_0} 2^{||\alpha||}\frac{{\theta-\gamma\choose \alpha}{\beta\choose 
\gamma}}{{\theta\choose \alpha}}=$$

$$\sum_{\gamma_1\in {\bf Z}_0}\cdots \sum_{\gamma_n\in {\bf Z}_0} \prod_{i=1}^n 2^{\alpha_i}\frac{{\theta_i-\gamma_i\choose \alpha_i}{\beta_i\choose 
\gamma_i}}{{\theta_i\choose \alpha_i}}=$$

$$ \prod_{i=1}^n \sum_{\gamma_i\in {\bf Z}_0}2^{\alpha_i}\frac{{\theta_i-\gamma_i\choose \alpha_i}{\beta_i\choose 
\gamma_i}}{{\theta_i\choose \alpha_i}}=$$
(by Corollary \ref{2})

$$ \prod_{i=1}^n \sum_{\gamma_i\in {\bf Z}_0}2^{\beta_i}\frac{{\theta_i-\gamma_i\choose \beta_i}{\alpha_i\choose 
\gamma_i}}{{\theta_i\choose \beta_i}}=$$

$$\sum_{\gamma\in {\bf Z}^n_0} 2^{||\beta||}\frac{{\theta-\gamma\choose \beta}{\alpha\choose 
\gamma}}{{\theta\choose \beta}}.$$

\begin{lm} \label{20} 
$$w(X)\ne 0\Rightarrow T(X)=0,$$
$$T(X)=T(\overline X).$$
\end{lm}

{\bf Proof.} It is enough to prove these statements for Weyl  monoms. Suppose that $X$ has a form 
$x^{\alpha}\der^{\beta}.$ If $\alpha\ne \beta,$ then
$T(X)=0.$ Therefore,  
$$T(X)=0=T(\overline{X}).$$
If $\alpha=\beta,$ then 
$$T(X)=\alpha!=\beta!=T(\overline{X}.$$

\begin{lm}\label{21} If $w(X)+w(Y)\ne 0,$ then $(X,Y)=0.$
\end{lm}

{\bf Proof.}  By Lemma \ref{20} 
$$w(X)+w(Y)\ne 0\Rightarrow T(X\circ Y)=0.$$
Therefore,
$$(X,Y)=T(X\circ Y)=0$$
if $w(X)+w(Y)\ne 0.$

\begin{lm} \label{4} Let $\alpha,\beta,\alpha',\beta'\in {\bf Z}^n.$ Then
$$(x^\alpha\der^\beta, x^{\alpha'}\der^{\beta'})=
\theta!\sum_{\gamma\in {\bf Z}_0^n}\frac{{\theta-\gamma\choose \alpha}{\beta\choose \gamma}}{{\theta\choose \alpha}},$$
if $\alpha'=\theta-\alpha$ and $\beta'=\theta-\beta$ for some $\theta\in {\bf Z}_0^n$ and 
$$(x^\alpha\der^\beta, x^{\alpha'}\der^{\beta'})=0$$
if such $\theta$ does not exist.  
\end{lm}

{\bf Proof.}  By Lemma \ref{21} it is enough to consider the case $w(X)+w(Y)=0.$
 
Suppose that  $X$ and $Y$ have forms $x^\alpha\der^\beta$ 
and $x^{\alpha'}\der^{\beta'}.$  Then $\alpha-\beta+\alpha'-\beta'=0.$ Set 
$\theta= \alpha+\alpha'.$ Then $\beta+\beta'=\theta$ and 
$$(X,Y)=T(x^{\alpha}\der^{\beta}\circ  x^{\alpha'}\der^{\beta'})=T(\sum_{\gamma\in {\bf Z}^n_0} {\beta\choose \gamma}{\alpha'\choose \gamma}\gamma! x^{\theta-\gamma}\der^{\theta-\gamma})=$$
$$  \sum_{\gamma\in {\bf Z}^n_0} {\beta\choose \gamma}{\alpha'\choose \gamma}\gamma! {(\theta-\gamma}!) =\theta!\sum_{\gamma\in {\bf Z}_0^n}\frac{{\theta-\gamma\choose \alpha}{\beta\choose \gamma}}{{\theta\choose \alpha}}.$$

\bigskip

{\bf Proof of Theorem \ref{Frobenius}}. Let us prove (\ref{51}).
 Let  $X\in A_{n,\omega}$ and $X'\in A_{n,\omega'}.$ 
By Lemma \ref{21}  $(X,Y)=0$ if $\omega'\ne -\omega.$ Therefore, in this case 
symmetric property of $(\;\,,\;)$ is evident.  Now assume that 
$\omega'=-\omega.$ This means that there exists some $\theta\in {\bf Z}^n$ such that 
$\theta=\alpha+\alpha'=\beta+\beta'.$ So, by Lemma \ref{4} and Lemma \ref{3} the symmetric property (\ref{51}) is established. 

Let us prove that  (\ref{51})  can be re-written in the  form (\ref{52}). 
If $l(X)+l(Y)\ne 0,$ then $w(X)+w(Y)\ne 0,$ and, by Lemma \ref{21} 
$$\sqrt{2}^{l(X)}(X,Y)=0=\sqrt{2}^{l(Y)}(Y,X).$$
If $l(X)+l(Y)=0,$ then 
$$(X,Y)=2^{l(Y)}(Y,X)\Rightarrow \sqrt{2}^{-l(Y)}(X,Y)=\sqrt{2}^{l(Y)}(Y,X)\Rightarrow 
\sqrt{2}^{l(X)}(X,Y)=\sqrt{2}^{l(Y)}(Y,X).$$
So, (\ref{52}) is established.
 
 By Lemma \ref{20} 
$$T(X\circ Y)=T(\overline{X\circ Y})=T(\overline{Y}\circ \overline{X}).$$
Therefore, (\ref{53}) is true. 
%$$(X,Y)= (\overline{Y},\overline{X}).$$

Associativity of bilinear form $(\;\,,\;)$ follows from associativity of Weyl algebra. 

Let $$I=\{ X\in A_n | (X,Y)=0, \quad  \forall Y\in A_n\}$$ 
be annulator of bilinear form $(\;\,,\;).$ We will prove that $I$ is two-sided ideal of $A_n.$ 

By associativity property for any  $X\in I$ and $Y, Y'\in A_n$ we have 
$$ (X\circ Y,Y')=(X, Y\circ Y')=0.$$
Thus, $X\circ Y\in I.$ So, $I$ is right-ideal of $A_n.$

Assume that  $X\in I$ is homogeneous, $X\in I\cap A_{n,l}$ for some  $l.$ 
Let $Y,Y'\in A_n$ are any homogeneous elements. 
By symmetric property 
$$(Y\circ X,Y')=2^{l(Y')-l(Y\circ X)} (Y',Y\circ X).$$
By associativity 
$$(Y',Y\circ X)=(Y'\circ Y,X).$$
By symmetric property 
$$(Y'\circ Y,X)=2^{-l(Y'\circ Y) 
+l(X)}(X,Y'\circ Y).$$
Since $X\in I,$
$$(X,Y'\circ Y)=0,$$
Hence
$$(Y\circ X,Y')=0$$
for any homogeneous elements  $X\in I,  \; Y,Y'\in A_n.$ Therefore by linearity
$$(Y,\circ X,Y')=0$$
for any $Y,Y'\in A_n$  and $X\in I.$ Thus, $Y\circ X\in I$ for any $Y\in A_n.$ 
This means that $I$ is left-ideal.

So, we have proved that $I$ is two-sided ideal of $A_n.$ 

Since $A_n$ is simple algebra, $I=0$ or $I=A_n.$ The second case is not 
possible: for example,
$$(x_1,\der_1)=T(x_1\der_1)=1\ne 0.$$
So,  $I=0.$ It  means that the form $( \; ,\,)$ is non-degenerate. 

\section{Euclid structure on Weyl algebra}

In this section we prove Theorem \ref{Euclid}. 

\begin{lm} \label{98}  Let 
$$\mu_{i,j}^{(a)}(t)=\sum_{i_1=0}^i {i\choose i_1}{a+j\choose a+i_1} t^{i_1}$$
and $M^{(a)}(t)=(\mu_{i,j}^{(a)}(t) )_{0\le i,j\le k}$ be $(k+1)\times (k+1)$-matrix. Then 
$$det\, M^{(a)}(t)=t^{{k+1\choose 2}}.$$
\end{lm}

{\bf Proof.} Let $M_1=(\nu_{i,j})_{0\le i,j \le k}$ be $(k+1)\times (k+1)$-matrix with $\nu_{i,j}=(-1)^{i+j}{i\choose j}.$ Then $M\cdot M_1=(\eta_{i,j})$ is $(k+1)\times (k+1)$-matrix such that 
$$\eta_{i,j}=\left\{\begin{array}{cc}
0& \mbox{ if $i>j$}\\
t^i&\mbox{ if $i=j\ge 0$}\end{array}\right.$$
Since
$$\nu_{i,j}=\left\{\begin{array}{cc}
1&\mbox{ if $i=j$}\\
0&\mbox{ if $i<j$}
\end{array}\right. ,
$$
it is clear that 
$$det\,M_1=1.$$
Therefore,
$$det\,M=det\,M\cdot M_1=t^{\sum_{i=0}^k i}=t^{{k+1\choose 2}}.$$

\begin{lm} \label{100} For any non-negative integers $a,b,c,$
$$\sum_{j=0}^a{a\choose j}{a+b-j\choose a+c}t^j=\sum_{i=0}^a{a\choose i}{b\choose i+c}(t+1)^i.$$
\end{lm}

{\bf Proof.} If $b<c$ both parts of this relation are $0.$ Consider the case $b\ge c.$ 
Imagine that    we have $a$ women and $b$ men and we like to construct  a group with $b-c$ persons such that $j$ women are selected as leaders.  We can do that in two ways. 

On the one hand we select $i$ women from  $a$ women, it can be done ${a\choose i}$ ways,
and among them select leaders, it can be done in ${i\choose j}$ ways, and select 
 $b-c-i$ men from $b$ men in ${b\choose b-c-i}$ ways. So, 
$\sum_{i=0}^a {i\choose j}{a\choose i}{b\choose b-c-i}$ is the number of ways to construct a group with $b-c$ persons with $j$ women leaders.

On the other hand we can select leaders first, it can be done in ${a\choose j}$ ways, and then select $b-c-j$ persons from $b+a-j$ persons, it can be done in ${a+b-j\choose a+c}$ ways. Hence, 
${a\choose j}{a+b-j\choose a+c}$ is the number of ways how to construct our group.

So, we establish that 
$$\sum_{i=0}^a {i\choose  j}{a\choose i}{b\choose b-c-i}={a\choose j}{a+b-j\choose a+c}.$$
for any $0\le j\le a.$ Therefore,
$$\sum_{j=0}^a {a\choose j}{a+b-j\choose a+c}t^j=
\sum_{j=0}^a \sum_{i=0}^a {i\choose  j}{a\choose i}{b\choose b-c-i}t^j=
\sum_{i=0}^a (t+1)^i {a\choose i}{b\choose b-c-i}.$$

\begin{lm} \label{101} Let 
$$\tilde d_{a,b}^{(c)}(x,y,z)=\sum_{i=0}^b {b\choose i}{a\choose i+c} x^{a-i}y^{b-i}(z+x y)^i,$$
$$d_{a,b}^{(c)}(x,y,z)=\sum_{i=0}^b {b\choose i}{a+b-i\choose b+c} x^{a-i}y^{b-i}z^i.$$
Then 
$$d_{a,b}^{(c)}(x,y,z)=\tilde d_{a,b}^{(c)}(x,y,z).$$
\end{lm}

{\bf Proof.} Take in Lemma \ref{100} $t=z/(xy)$ and multiply both paths of this relation to 
$x^a y^b.$ We obtain
$$\sum_{j=0}^a {a\choose j}{a+b-j\choose a+c}x^{a-j}y^{b-j}z^j=
\sum_{i=0}^a{a\choose i}{b\choose i+c}x^{a-i}y^{b-i}(xy+z)^i.$$
Change here $a$ to $b$ and $b$ to $a.$ We  obtain that $d_{a,b}^{(c)}(x,y,z)=\tilde d_{a,b}^{(c)}(x,y,z).$

\begin{lm}\label{102}  Let $\tilde M^{(a,k)}(x,y,z)$ be $(k+1)\times (k+1)$-matrix with $(i,j)$-components 
$\tilde d^{(a)}_{a+j,i}(x,y,z),$ where $0\le i,j\le k.$ Then 
$$det\,\tilde M^{(a,k)}(x,y,z)=x^{a(k+1)} (xy+z)^{{k+1\choose 2}}.$$
\end{lm}

{\bf Proof.} Note that 
$$\tilde d_{a+j,i}^{(a)}(x,y,z)=\sum_{i_1=0}^{i}{i\choose i_1}{a+j\choose a+i_1}x^{a+j-i_1}y^{i-i_1}(z+xy)^{i_1}
=$$
$$x^{a+j}y^i\sum_{i_1=0}^{i}{i\choose i_1}{a+j\choose a+i_1}(1+z/(xy))^{i_1}=\mu_{i,j}^{(a)}(1+z/(xy)).$$
(definition of $\mu_{i,j}^{(a)}$  is given in Lemma \ref{98}). 
Therefore, by Lemma \ref{98}
$$det\,\tilde M^{(a,k)}(x,y,z)=\prod_{j=0}^k x^{a+j}\prod_{i=0}^k y^i (1+z/(xy))^{{k+1\choose 2}}=x^{a(k+1)} (xy+z)^{{k+1\choose 2}}.$$

\begin{lm}\label{103} Let 
$$\eta_{i,j}^{(a,b)}=\sum_{i_1=0}^{b+i} {b+i\choose i_1}{b+j\choose i_1}i_1!(a+b+i+j-i_1)!.$$
Let $N^{(a,k)}$ be $(k+1)\times (k+1)$-matrix with $(i,j)$-components  $\eta_{i,j}^{(a,0)}.$ Then
$$det\,N^{(a,k)}>0.$$
\end{lm}

{\bf Proof.} 
Note that 
$$\eta_{i,j}^{(a,b)}=(a+i)!(b+j)!\sum_{i_1=0}^{b+i} {b+i\choose i_1}{a+b+i+j-i_1\choose a+i}.$$
Therefore,
$$\eta_{i,j}^{(a,0)}=(a+i)!\,j!\,\sum_{i_1=0}^{i} {i\choose i_1}{a+i+j-i_1\choose a+i}=
(a+i)!\, j!\, d_{a+j,i}^{(a)}(1,1,1).$$

Let $M^{(a,k)}(x,y,z)$ be $(k+1)\times (k+1)$-matrix with $(i,j)$-components 
$d_{a+j,i}^{(a)}(x,y,z),$ where $0\le i,j\le k.$ Then by Lemma \ref{101}
$$det\,M^{(a,k)}(x,y,z)=det\,\tilde M^{(a,k)}(x,y,z).$$
Therefore by Lemma \ref{102} 
$$det\,N^{(a,k)}=\prod_{i,j=0}^k i!\,(a+j)! \, det\,M^{(a,k)}(1,1,1)= \prod_{i,j=0}^k i!\,(a+j)! \, 2^{{k+1\choose 2}}>0$$

\begin{lm} \label{104}
Let $X=\sum_{i\ge 0} \lambda_i X_i\in A_1,$ $X_i=x^{a+i}\der^i,$ where number of non-zero coefficients $\lambda_i$ is finite and $a\ge 0.$Then $\langle X,X\rangle\ge 0$ and the condition $\langle X,X\rangle=0$ implies that $X=0.$
\end{lm}

{\bf Proof.}  We have 
$$\langle X,X\rangle=\sum_{i,j\ge 0} \lambda_i\lambda_j \langle X_i,X_j\rangle=
$$
$$\sum_{i,j\ge 0}\lambda_i\lambda_j \sum_{i_1=0}^{i}{i\choose i_1}{j\choose i_1}i_1! (a+i+j-i_1)!=
$$
$$\sum_{i,j\ge 0} \lambda_i\lambda_j\eta_{i,j}^{(a,0)}.$$
Therefore, by Sylvester's criterion our Lemma is equivalent to the condition 
$$det\, N^{(a,k)}>0$$
for all $k\ge 0.$ By Lemma \ref{104} it is true. 

\bigskip

{\bf Proof of Theorem \ref{Euclid}}.  By (\ref{53})
$$\sqrt{2}^{l(X)}(X,\overline{Y})=\sqrt{2}^{l(X)}(Y,\overline{X}).$$
Therefore, if $l(X)=l(Y),$ then 
$$\langle X,Y\rangle =\langle Y,X\rangle.$$
If $l(X)\ne l(Y),$ then $w(X)\ne w(Y),$ and by Lemma \ref{21}
$$\langle X,Y\rangle =0=\langle Y,X\rangle.$$
So, (\ref{61}) is true. 

Let us prove (\ref{62}).  By (\ref{53}) 
$$(\overline{X},Y)=(\overline{Y},X).$$
By (\ref{51})
$$(\overline{Y},X)=2^{l(X)}(X,\overline{Y}).$$
Therefore, 
$$2^{l(X)}(X,\overline{Y})=(\overline{X},Y),$$
and,
$$\sqrt{2}^{l(X)}(X,\overline{Y})=
\sqrt{2}^{\; -l(X)}(\overline{X},Y)
=\sqrt{2}^{l(\overline{X})}(\overline{X},Y).
$$
We obtain (\ref{62}).

Let us prove (\ref{63}). 
If $l(a)+l(X)\ne l(Y),$ then $w(a)+w(X)\ne w(Y),$ and 
$w(X)\ne w(\overline{a})+\omega(Y).$ Therefore by Lemma \ref{21}
$$\langle a\circ X,Y\rangle=0=\langle X,\overline a\circ Y\rangle.$$
if $l(a)+l(X)\ne l(Y).$

Suppose now that $l(a)+l(X)= l(Y).$ Then
\begin{equation}\label{mmm}
l(a\circ X)/2+l(X)/2+l(\overline{Y})=l(a)/2+l(X)/2+l(X)/2-l(Y)=-l(a)/2.
\end{equation}
By (\ref{51})
$$\langle a\circ X,Y\rangle=\sqrt{2}^{ l(a\circ X)} (a\circ X,\overline{Y})=
\sqrt{2}^{ l(a\circ X)} 2^{l(\overline{Y})}(\overline{Y},a\circ X)=$$
By associative property of Frobenius form 
$$(\overline{Y},a\circ X)=(\overline{Y}\circ a, X).$$
By (\ref{51})
$$(\overline{Y}\circ a, X)=2^{l(X)}(X,\overline{Y}\circ a).$$
Therefore,  
$$\langle a\circ X,Y\rangle=2^{l(a\circ X)/2+l(X)+l(\overline{Y})}(X,\overline{Y}\circ a).$$
So, by (\ref{mmm}), 
$$\langle a\circ X,Y\rangle=2^{\; -l(a)/2} \sqrt{2}^{l(X)}(X,\overline{\overline{a}\circ Y})
=\sqrt{2}^{\; -l(a)}\langle X,\overline{a}\circ Y\rangle.$$
Relation (\ref{63}) is proved. 

Proof of  (\ref{64}) is similar to the proof of (\ref{63}).

Recall that by  $X_{\omega}$ we denote a projection of $X\in A_n$ on $A_{n,\omega}.$ By Lemma \ref{21}
$$\omega\ne \omega'\Rightarrow \langle X_{\omega}, Y_{\omega'}\rangle=0.$$
Therefore
$$\langle X,Y\rangle =\sum_{\omega\in {\bf Z}^n} 
\langle X_\omega,Y_\omega\rangle.$$

Recall that
$$A_1^{(i)}=\{ x_i^{\alpha_i}\der_i^{\beta_i} | \alpha_i,\beta_i\in {\bf Z} \}$$
is a  Weyl algebra in one variable $x_i.$
Suppose that  $X$ has a form $x^\alpha \der^\beta.$
Let 
$$X^{(i)}=x_i^{\alpha_i} \der_i^{\beta_i}.$$
Then 
$$X=X^{(1)}\cdots X^{(n)}$$
and the map $X\mapsto X^{(1)}\otimes \cdots X^{(n)}$  gives us an isomorphism of associative algebras  
$$A_n\cong A_1^{(1)}\otimes \cdots A_1^{(n)}.$$

Note that 
$$\langle X,X\rangle =\langle X^{(1)}, X^{(1)}\rangle \cdots \langle X^{(n)}, X^{(n)}\rangle.$$
If 
$$\langle X^{(i)}, X^{(i)}\rangle \ge 0, \qquad  i=1,\ldots,n, $$
then 
$$\langle X,X\rangle \ge 0.$$
If the condition 
$$\langle X^{(i)},X^{(i)}\rangle =0$$
implies that $X^{(i)}=0,$ then the condition $\langle X, X\rangle =0$ implies that 
$\langle X^{(i)},X^{(i)}\rangle=0$ for some $i,$ and, $X=X^{(1)}\cdots X^{(n)}=0.$

Therefore, to establish the positive definity of the bilinear form $\langle\;,\;\rangle$ on $A_n$ it is enough to prove it for the case $n=1.$ Moreover we can assume that $X\in A_1$ is homogeneous.
So, we have to prove that $\langle X,X \rangle \ge 0$ for 
$$X=\sum_{i\in {\bf Z}} \lambda_i x^{a+i}\der^{b+i} $$
where $a,b\in {\bf Z}_0,$ such that $a+i\ge 0, b+i\ge 0.$ 
By (\ref{63}) and (\ref{64})  it is enough to prove it  for $a\ge 0$ and $b=0.$ 
By Lemma \ref{104}
$$X=\sum_{i\ge 0} \lambda_i x^{a+i}\der^i, \; a\ge 0\Rightarrow \langle X,X\rangle \ge 0,$$
and 
$$\langle X,X\rangle =0, X=\sum_{i\ge 0} \lambda_i x^{a+i}\der^i, \; a\ge 0\Rightarrow \lambda_i=0, i\ge 0.$$
So,  positive definity of the scalar product $\langle X,X\rangle$ is established.

\end{document}